\documentclass[12pt,leqno]{article}
\usepackage{amsfonts}
\usepackage{amsmath, amsthm, amsfonts, amssymb, color}
\usepackage{mathrsfs}
\usepackage{color}
\theoremstyle{definition}

\newcommand{\scr}[1]{\mathscr #1}
\definecolor{wco}{rgb}{0.5,0.2,0.3}

\numberwithin{equation}{section} \theoremstyle{remark}

\newcommand{\ua}{\uparrow}

\title{{\bf  Extrinsic Derivative Formula for Distribution Dependent SDEs }\footnote{The author is supported by  NNSFC (12301180) and Research Centre for Nonlinear Analysis at Hong Kong PolyU.} }
\author{
{\bf Panpan Ren}\\
Mathematics Department,  City University of Hong Kong,  Hong Kong,  China\\
panparen@cityu.edu.hk
 }
\begin{document}
\allowdisplaybreaks
\def\R{\mathbb R}  \def\ff{\frac} \def\ss{\sqrt} \def\B{\mathbf
B}
\def\N{\mathbb N} \def\kk{\kappa} \def\m{{\bf m}}
\def\ee{\varepsilon}\def\ddd{D^*}
\def\dd{\delta} \def\DD{\Delta} \def\vv{\varepsilon} \def\rr{\rho}
\def\<{\langle} \def\>{\rangle}
  \def\nn{\nabla} \def\pp{\partial} \def\E{\mathbb E}
\def\d{\text{\rm{d}}} \def\bb{\beta} \def\aa{\alpha} \def\D{\scr D}
  \def\si{\sigma} \def\ess{\text{\rm{ess}}}\def\s{{\bf s}}
\def\beg{\begin} \def\beq{\begin{equation}}  \def\F{\scr F}
\def\Ric{\mathcal Ric} \def\Hess{\text{\rm{Hess}}}
\def\e{\text{\rm{e}}} \def\ua{\underline a} \def\OO{\Omega}  \def\oo{\omega}
 \def\tt{\tilde}\def\[{\lfloor} \def\]{\rfloor}
\def\cut{\text{\rm{cut}}} \def\P{\mathbb P} \def\ifn{I_n(f^{\bigotimes n})}
\def\C{\scr C}      \def\aaa{\mathbf{r}}     \def\r{r}
\def\gap{\text{\rm{gap}}} \def\prr{\pi_{{\bf m},\varrho}}  \def\r{\mathbf r}
\def\Z{\mathbb Z} \def\vrr{\varrho} \def\ll{\lambda}
\def\L{\scr L}\def\Tt{\tt} \def\TT{\tt}\def\II{\mathbb I}
\def\i{{\rm in}}\def\Sect{{\rm Sect}}  \def\H{\mathbb H}
\def\M{\mathbb M}\def\Q{\mathbb Q} \def\texto{\text{o}} \def\LL{\Lambda}
\def\Rank{{\rm Rank}} \def\B{\scr B} \def\i{{\rm i}} \def\HR{\hat{\R}^d}
\def\to{\rightarrow}\def\l{\ell}\def\iint{\int}\def\gg{\gamma}
\def\EE{\scr E} \def\W{\mathbb W}
\def\A{\scr A} \def\Lip{{\rm Lip}}\def\S{\mathbb S}
\def\BB{\scr B}\def\Ent{{\rm Ent}} \def\i{{\rm i}}\def\itparallel{{\it\parallel}}
\def\g{{\mathbf g}}\def\Sect{{\mathcal Sec}}\def\T{\mathcal T}\def\BB{{\bf B}}
\def\f{\mathbf f} \def\g{\mathbf g}\def\BL{{\bf L}}  \def\BG{{\mathbb G}}
\def\Bd{{D^E}} \def\BdP{D^E_\phi} \def\Bdd{{\bf \dd}} \def\Bs{{\bf s}} \def\GA{\scr A}
\def\Bg{{\bf g}}  \def\Bdd{\psi_B} \def\supp{{\rm supp}}\def\div{{\rm div}}
\def\ddiv{{\rm div}}\def\osc{{\bf osc}}\def\1{{\bf 1}}\def\BD{\mathbb D}\def\GG{\Gamma}
\def\H{{\bf H}}
\maketitle

\begin{abstract} A Bismut type formula is established for the extrinsic derivative  of distribution dependent SDEs (DDSDEs). The main result is illustrated by  non-degenerate DDSDEs  with space-time singular drift, as well as   degenerate  DDSDEs with weakly   monotone coefficients.
  \end{abstract} \noindent


 \vskip 2cm

 \section{Introduction}

 Let $k\in [0,\infty)$ and $\scr P_k$ be the set of all probability measures $\mu$ on $\R^d$ having finite $k^{th}$-moment; i.e.
 $$\mu(|\cdot|^k):=\int_{\R^d}|x|^k\mu(\d x)<\infty.$$ In particular, $\scr P_0=\scr P$ is the set of all probability measures on $\R^d$.
 It is well known that $\scr P_k$ is a complete metric space under the weighted total variance  distance
 $$\|\mu-\nu\|_{k,var}:=\sup_{|f|\le 1+|\cdot|^k} |\mu(f)-\nu(f)|,$$ see for instance \cite{V}.

  Consider the following distribution dependent SDE (DDSDE) on $\R^d$:
 \beq\label{E0} \d X_t= b_t(X_t,\L_{X_t})\d t+\si_t(X_t)\d W_t,\ \ t\in [0,T],\end{equation}
 where $T>0$ is a fixed time, $\L_{X_t}$ is the distribution of $X_t$, $W_t$ is the $d$-dimensional Brownian motion on a probability base $(\OO,\{\F_t\}_{t\in [0,T]},\F,\P)$, and
 $$b: [0,T]\times \R^d\times \scr P_k\to \R^d,\ \ \si: [0,T]\times \R^d\to \R^d\otimes \R^d$$
 are measurable. We call \eqref{E0} well-posed for distributions in $\scr P_k$, if for any $\F_0$-measurable initial value $X_0$ with $\E|X_0|^k<\infty$ (resp. any initial distribution $\mu\in \scr P_k$),
the SDE has a unique strong (resp. weak) solution with $\L_{X_\cdot}\in C([0,T]; \scr P_k).$
In this case, to emphasize the initial distribution we denote the solution by $X_t^\mu$ if $\L_{X_0}=\mu,$
and denote $P_t^*\mu=\L_{X_t^\mu}.$

For any $f\in \B_b(\R^d)$, let
$$P_tf(\mu)=\E[f(X_t^\mu)],\ \ \mu\in \scr P_k.$$
 To characterize the regularity of non-linear Fokker-Planck equations, following type Bismut   formula  have been established in \cite{RW19, BRW, HW21, W23} for the intrinsic/Lions derivative of $\mu\mapsto P_tf(\mu)$:
 $$DP_tf(\mu)=\E\big[f(X_t^\mu) M_t],\ \ f\in \B_b(\R^d),$$ where $M_t$ is a martingale explicitly given by the solution of \eqref{E0}. The crucial point of this type formula is that the derivative of $P_tf$ is represented by $f$ rather than the derivative of $f$, so that it implies regularity estimates for the distribution of solutions.   In the study of particle systems, the intrinsic derivative describes the motion while the extrinsic derivative describes the birth-death of particles, see for instance \cite{KLV,RW20} and references therein. 

 In this paper, we aim to establish the Bismut type  derivative formula  the same type formula for the  extrinsic derivative of $P_tf$,  which is not yet available so far.

\beg{defn} Let $f$ be a continuous function on $\scr P_k$. 
 \beg{enumerate} \item[(1)] We call $f$ extrinsically  differentiable, if for any $\mu\in \scr P_k$, the convex extrinsic derivative
$$\tt D^E_x f(\mu):=\lim_{\vv \downarrow 0}\ff{f((1-\vv)\mu+\vv\dd_x)-f(\mu)}\vv\in\R,\ \ x\in\R^d $$
exists, where $\dd_x$ is the Dirac measure at point $x$.
\item[(2)]  We denote $f\in C^{E,1}(\scr P_k)$ if $f$ is extrinsically differentiable, and
$$(x,\mu)\in \R^d\times \scr P_k \mapsto \tt D^E  f(\mu)(x):= \tt D^E_x f(\mu)$$ is  continuous.
\item[(3)] We write $f\in C_K^{E,1}(\scr P_k)$, if $f\in C^{E,1}(\scr P_k)$ and for any compact set
$\scr K\subset \scr P_k$ there exists a constant $c>0$ such that
$$\sup_{\mu\in \scr K}|\tt D^E f(\mu)(x)|\le c(1+|x|^k),\ \ x\in\R^d.$$\end{enumerate}
\end{defn}

As shown in \cite[Lemma 3.2]{RW21}, when $f\in C_K^{E,1}$ the function $x\mapsto \tt D^E_x f(\mu)$
becomes the linear functional derivative, i.e. for any $\nu\in  \scr P_k$,
\beq\label{LF} f(\mu)-f(\nu)=\int_0^1 (\mu-\nu)\big(\tt D^E f(r\mu+(1-r)\nu)\big)\,\d r.\end{equation}

To calculate $\tt D^E P_tf$, we will use the semigroup of the decoupled SDE associated with \eqref{E0}:
\beq\label{DC} \d X_t^{\mu,x}= b_t(X_t^{\mu,x},P_t^*\mu)\d t+ \si_t(X_t^{\mu,x})\d W_t,\ \ X_0^{\mu,x}=x, t\in [0,T],\end{equation} where $x\in \R^d$ and $\mu\in \scr P_k.$ Let $P_t^\mu$ be the associated semigroup, i.e.
$$P_t^\mu f(x):=\E[f(X_t^{\mu,x})],\ \ t\in [0,T], x\in\R^d.$$
In general, for a probability measure $\nu$ on $\R^d$, let
$$P_t^\mu f(\nu):=\int_{\R^d} P_t^\mu f(x)\nu(\d x),\ \ (P_t^\mu)^*\nu:=\int_{\R^d}\L_{X_t^{\mu,x}}\nu(\d x).$$

We make the following assumption.

\beg{enumerate}\item[{\bf (H)}] Let $k\in [0,\infty)$. The following two conditions hold.
\item[$(H_1)$] For any $\gg_\cdot\in C([0,T];\scr P_k)$, the SDE
\beq\label{DC'} \d X_t^{\gg_\cdot,x}= b_t(X_t^{\gg_\cdot,x},\gg_t)\d t+ \si_t(X_t^{\gg_\cdot,x})\d W_t,\ \ X_0^{\gg_\cdot,x}=x, t\in [0,T],\end{equation}  is well-posed, and  there exist constants $c >0$ and $p\geq k$ independent of $\gg_\cdot$ such that
\beq\label{DC1} \E\big[|X_t^{\gg_\cdot,x}|^{p}\big]\le c \big(1+|x|^{p}+\int_0^t\gg_s(|\cdot|^k)^{\ff p k}\d s\big),\ \ t\in [0, T], \\ x\in \R^d,\end{equation} where 
$\gg_t(|\cdot|^k)^{\ff p k}:=1$ when $k=0.$
\item[$(H_2)$] The drift  $b_t(x,\mu)$ is formulated as
$$b_t(x,\mu)= b_t^{(0)}(x)+ \si_t(x) b_t^{(1)}(x,\mu),$$
where
 $$b^{(0)}: [0,T]\times \R^d\to \R^d,\ \ b^{(1)}: [0,T]\times \R^d\times\scr P_k\to \R^d$$ are measurable, $b_t^{(1)}(x,\cdot)\in C^{E,1}(\scr P_k),$          there exists   $K\in L^2([0,T];(0,\infty))$ such that
\beq\label{DC2}   \inf_{c\in \R^d}|\tt D^E b_t^{(1)}(x,\mu)(y)-c|\le K_t(1+ |y|^k),\ \ t\in [0,T],\  x,y\in \R^d,\ \mu \in \scr P_k, \end{equation} and there exists an increasing function $\aa: (0,\infty)\to (0,\infty)$ with   $\aa(\vv)\to 0$ as $\vv\to 0$ such that
\beq\label{DC3} \beg{split}&|\tt D^E b_t^{(1)}(x,\mu)(y)- \tt D^E b_t^{(1)}(x,\nu)(y)|\\
&\le K_t\aa(\|\mu-\nu\|_{k,var})\big(1+|y|^k+\mu(|\cdot|^k)+\nu(|\cdot|^k)\big),\ \ x,y\in \R^d,\ \mu,\nu\in \scr P_k.\end{split}\end{equation} \end{enumerate}

The condition $(H_1)$ can be checked easily by existing results concerning the well-posedness and moment estimates. So, in applications  the key point is to verify $(H_2)$ for the distribution dependence of $b$. A simple example satisfying $(H_2)$ is
$$b_t^{(1)}(x,\mu)=F_t(x, \mu(h)),$$
where for some $n\ge 1$,
$$F: [0,T]\times \R^d\times \R^n\to\R^d,\ \ h:\R^d\to\R^n$$ are measurable   such that for any $(t,x)\in [0,T]\times \R^d$,
$$\|\nn  F_t(x,\cdot)\|\lor  \|\nn^2  F_t(x,\cdot)\| \le C,\ \ \ |h(x)| \le C(1+|x|^k)$$ hold for some constant $C>1$. In this case, we have
\beg{align*}&\inf_{c\in\R}|\tt D^E b_t^{(1)}(x,\mu)(y)-c|\le \big|\nn F_t(x,\cdot)(\mu(h)) h(y)\big|  \le C(1+ |y|^k),\\
&|\tt D^E b_t^{(1)}(x,\mu)(y)- \tt D^E b_t^{(1)}(x,\nu)(y)|\\
&= \big|[\nn F_t(x,\cdot)(\mu(h))]\big(h(y)-\mu(h)\big)- [\nn F_t(x,\cdot)(\nu(h))]\big(h(y)-\nu(h)\big)\big|\\
&\le \big\|\nn^2 F_t(x,\cdot)\big\|_\infty \big|\big(\mu(h)-\nu(h)\big)\big| \cdot \big|\big(h(y)-\mu(h)\big)\big|+ \big|[\nn F_t(x,\cdot)(\nu(h))](\mu(h)-\nu(h))\big|\\
&\le C^2 \big(2+\mu(|\cdot|^k)+|y|^k\big)\|\mu-\nu\|_{k,var}+C\|\mu-\nu\|_{k,var}.\end{align*}
So, $(H_2)$ holds.

\

Let $\D_k$ be the class of measurable functions $f$ on $\R^d$ such that 
$|f|\le c(1+|\cdot|^k)$ for some constant $c>0$.
The main result of this paper is the following.

\beg{thm}\label{T1} Assume {\bf (H)}. Then the DDSDE $\eqref{E0}$ is well-posed for distributions in $\scr P_k$, and the following assertions hold.
\beg{enumerate} \item[$(1)$] For any $\mu,\nu\in \scr P_k$, the equation
\beq\label{ET}\beg{split} &\eta_t^{\mu,\nu}= P_t^\mu\big(\tt D^E b_t(X_t^\mu, P_t^*\mu)\big)(\nu)- P_t\big(\tt D^E b_t(X_t^\mu, P_t^*\mu)\big)(\mu) \\
&\quad \qquad + \E\bigg[\tt D^E b_t(z, P_t^*\mu)(X_t^\mu)
\int_0^t \big\<\eta_s^{\mu,\nu}, \d W_s\big\>\bigg]_{z=X_t^\mu},\ \ t\in [0,T]\end{split}\end{equation}
has a unique solution satisfying
\beq\label{ET2}  \|\eta_t^{\mu,\nu}\|_\infty \le cK_t\big(1+ (\mu+\nu)(|\cdot|^{\theta})\big)
\e^{c\mu(|\cdot|^{k})^2},\ \ \mu,\nu\in \scr P_k,\ t\in [0,T]\end{equation}
for some constant $c>0.$
\item[$(2)$] For any $t\in (0,T]$ and $f\in \D_{k},$
\beq\label{BS} \tt D^E_\nu P_t f(\mu)= \int_{\R^d} \big(P_t^\mu f\big)\d (\nu-\mu)
+\E\bigg[f(X_t^\mu)\int_0^t\big\< \eta_s^{\mu,\nu},\ \d W_s\big\>\bigg].
\end{equation} Consequently,    there exist a  constant  $c >0$   such that
\beq\label{DE}\beg{split} & \sup_{|f|\le 1+|\cdot|^k}|\tt D^E_\nu P_t f(\mu)| 
\le c \|\mu-\nu\|_{k,var}\\
&\quad +   (1+\mu(|\cdot|^k))
\big(1+(\mu+\nu)(|\cdot|^k)\big) \e^{c\mu(|\cdot|^{k})^2}\ss t,
\ \ t\in [0,T].\end{split}\end{equation}
\end{enumerate}
\end{thm}

\paragraph{Remark 1.1.} For $\mu,\nu\in \scr P_k$, let
$$\tt D_\nu^E P_t^*\mu:={\D_k}\text{-}\lim_{\vv\downarrow 0} \ff{P_t^*((1-\vv)\mu+\vv \nu)-P_t^*\nu}\vv,$$
where the limit is defined with respect to $\D_k$: a sequence of signed measures $\{\phi_n\}_{n\ge 1}$ converges to $\phi$ with respect to $\D_k$,  if $(|\phi_n|+|\phi|)(|\cdot|^k)<\infty$ and
$$\lim_{n\to\infty} \phi_n(f)=\phi(f),\ \  f\in \D_k.$$
By \eqref{BS}, $\tt D_\nu^E P_t^*\mu$ exists and
$$\tt D_\nu^E P_t^*\mu=(P_t^\mu)^*(\nu-\mu)+ Z_t(y)(P_t^*\mu)(\d y),$$
where $Z: [0,T]\times \R^d\to\R^d$ is measurable such that
$$Z_t(X_t^\mu)=\E\bigg[\int_0^t \big\<\eta_s^{\mu,\nu}, \d W_s\big\>\bigg|X_t^\mu\bigg].$$

\

In Section 2, we present some lemmas, which will be used in Section 3 to prove Theorem \ref{T1}.
In Sections 4, we apply the main result to non-degenerate singular DDSDEs and degenerate  
models.

 \section{Lemmas}

 \beg{lem}\label{L0} Assume {\bf (H)}. Then the SDE $\eqref{E0}$ is well-posed for distributions in $\scr P_k$, and there exists a constant $C>0$ such that for any $\mu\in \scr P_k$,
 \beq\label{DC1*}    \E\big[|X_t^\mu|^p\big|\F_0\big]\le C \big(1+ |X_0^\mu|^p+\mu(|\cdot|^k)^{\ff p k}\big).\end{equation}
 \end{lem}

 \beg{proof} The proof is similar to that of \cite[Theorem 3.2]{W23b}  using the fixed point theorem.
By $(H_1)$, for any
 $$\gg_\cdot\in \C^\mu:=\big\{\gg_\cdot\in C([0,T];\scr P_k):\ \gg_0=\mu\big\},$$
 the SDE
\beq\label{EE}\d X_t^{\gg_\cdot}= b_t(X_t^{\gg_\cdot}, \gg_t)\d t +\si_t(X_t^{\gg_\cdot})\d W_t,\ \ t\in [0,T], X_0^{\gg_\cdot}=X_0^\mu\end{equation}
 has a unique solution, such that
 \beq\label{PP1} \E[|X_t^{\gg_\cdot}|^p|\F_0]\le c  (1+ |X_0^{\mu}|^p)+c \int_0^t
 \gg_s(|\cdot|^k)^{\ff p k}\d s,\ \ t\in [0,T].\end{equation}
So,  the map $\Phi: \C^\mu\to\C^\mu$ defined by
 $$\Phi_t(\gg):=\L_{X_t^{\gg_\cdot}},\ \ t\in [0,T]$$ satisfies
 \beq\label{PP2}  \E[|X_t^{\gg_\cdot}|^k] \le c^{\ff k p}  (1+ \mu(|\cdot|^k)) + \bigg(c\int_0^t \gg_s(|\cdot|^k)^{\ff p k}\d s\bigg)^{\ff k p},\ \ t\in [0,T].\end{equation}
 Then there exists a  constant $N_0\ge 1$, such that for any $N\ge N_0$, the set
 $$\C^\mu_N:=\Big\{\gg_\cdot\in C([0,T];\scr P_k): \sup_{t\in [0,T]}\e^{-Nt}\gg_t(|\cdot|^k)\le N(1+\mu(|\cdot|^k)\Big\}$$ is invariant under the map $\Phi$.
 Next, $(H_2)$ implies that $b_t(x,\cdot)\in C_K^{E,1}(\scr P_k)$. By \eqref{LF} we obtain 
\beq\label{CP0}\beg{split} &\big|b_t^{(1)}(x,\mu)-b_t^{(1)}(x,\nu)\big|
=\bigg|\int_0^1 \Big(\ff{\d}{\d r} b_t\big(x, r\mu+(1-r)\nu\big)\Big)\d r\bigg|\\
&=\bigg|\int_0^1\d r \int_{\R^d} \tt D^E b_t^{(1)}\big(x,r\mu+(1-r)\nu\big)(z) (\mu-\nu)(\d z)\bigg|\\
&\le K_t\|\mu-\nu\|_{k, var},\ \ t\in [0,T], \ \mu,\nu\in \scr P_k. \end{split}\end{equation}
Then for any $\gg_\cdot,\tt\gg_\cdot\in \C_N^\mu$, we reformulate the SDE \eqref{EE} as
$$\d X_t^{\gg_\cdot}= b_t(X_t^{\gg_\cdot}, \tt\gg_t)\d t +\si_t(X_t^{\gg_\cdot})\d \tt W_t,\ \ t\in [0,T], X_0^{\gg_\cdot}=X_0^\mu,$$
where by Girsanov's theorem,
$$\tt W_t:=W_t-\int_0^t \eta_s \d s,\ \ t\in [0,T]$$ is a Brownian motion under the weighted probability measure $R_T\d\P$, for $R_T$ defined as
\beg{align*} &\eta_s:=\big[b_s(X_s^{\gg_\cdot}, \tt\gg_s)-b_s(X_s^{\gg_\cdot},  \gg_s)\big],\\
&R_t:= \e^{\int_0^t\<\eta_s,\d W_s\>-\ff 1 2\int_0^t|\eta_s|^2\d s},\ \ s,t\in [0,T].\end{align*}
By the weak uniqueness, and H\"{o}lder's inequality,  we find a constant $c_1>0$ such that
\beg{align*}&\|\Phi_t(\gg)-\Phi_t(\tt\gg)\|_{k,var}=\sup_{|f|\le 1+|\cdot|^k} \big|\E[f(X_t^{\gg_\cdot})(R_t-1)]\big|\\
&\le \E[(1+|X_t^{\gg_\cdot}|^k)|R_t-1|]\le c_1 \E\Big[\Big(\E\big((1+|X_t^{\gg_\cdot}|^p)\big|\F_0\big)\Big)^{\ff k p}
\Big(\E\big(|R_t-1|^{\ff p{p-k}}\big|\F_0\big)\Big)^{\ff{p-k}p}\Big].\end{align*}
By \eqref{CP0} and $\gg_\cdot,\tt\gg_\cdot\in \C_N^\mu$, and noting that
$$|\e^r-1|\le  (\e^r+1)|r|,\ \ r\in\R,$$ we find a constant $c_2(N)>0$ such that
$$\Big(\E\big(|R_t-1|^{\ff p{p-k}}\big|\F_0\big)\Big)^{\ff{p-k}p} \le c_2(N)\bigg(\int_0^t K_s^2 \|\gg_s-\tt\gg_s\|_{k,var}^2\d s\bigg)^{\ff 1 2}.$$ Combining this with  \eqref{PP1} we arrive at
$$\|\Phi_t(\gg)-\Phi_t(\tt\gg)\|_{k,var}\le c_3(N) \bigg(\int_0^tK_s^2 \|\gg_s-\tt\gg_s\|_{k,var}^2\d s\bigg)^{\ff 1 2}.$$
Therefore, when $\ll>0$ is large enough the map
$\Phi$ is contractive on $\C_\mu^N$ under the complete metric
$$\rr_\ll(\gg_\cdot,\tt\gg_\cdot):=\sup_{t\in [0,T]}\e^{-\ll t}\|\gg_t-\tt\gg_t\|_{k,var},$$
and hence the SDE \eqref{E0} is well-posed for distributions in $\scr P_k$.

Finally, by \eqref{PP1} for the fixed point $\gg_t:=P_t^*\mu$, we obtain
$$\E[|X_t^{\mu}|^p|\F_0]\le c  (1+ |X_0^{\mu}|^p)+c \int_0^t
 \big(\E[|X_s^\mu|^k]\big)^{\ff p k}\d s,\ \ t\in [0,T].$$
 This implies   \eqref{DC1}  for some constant $C>0$.

 \end{proof}

Noting that $X_t^\mu$ solves \eqref{DC}
for initial value $X_0^\mu$ replacing $x$, we have
\beq\label{D**} P_t^*\mu:=\L_{X_t^\mu}=(P_t^\mu)^*\mu,\end{equation}
and \eqref{DC1} implies
\beq\label{DC1*} \E\bigg[\sup_{t\in [0,T]}|X_t^\mu|^p\Big|\F_0\bigg]= \E\bigg[\sup_{t\in [0,T]}|X_t^{\mu,x}|^p\Big|\F_0\bigg]_{x=X_0^\mu}\le c_p(1+|X_0^{\mu}|^p), \ \ \ \ p\in [1,\infty), \mu\in \scr P_k,\end{equation}
\beq\label{DC1**} \int_{\R^d}|\cdot|^p\d  (P_t^\mu)^*\nu\le c_p (1+\nu(|\cdot|^p)),\ \ p\in [1,k], \mu,\nu\in \scr P_k.
\end{equation}

\beg{lem}\label{L1} Assume {\bf (H)}. Simply denote $\pi_\vv= (1-\vv)\mu+\vv\nu, \vv\in (0,1).$  Then   there exists a constant $c(\mu,\nu)>0$ increasing in $(\mu+\nu)(|\cdot|^k) $  such that
$$\|(P_t^{\pi_\vv})^*\gg- (P_t^\mu)^*\gg\|_{k,var}\le  c(\mu,\nu)(1+\gg(|\cdot|^k))\vv,
  \ \ \vv\in (0,1), \gg\in \scr P_k. $$
\end{lem}

\beg{proof} Let $X_0$ be $\F_0$-measurable with $\L_{X_0}=\gg$. Consider the SDE
\beq\label{E01}\d X_t= \big\{b_t^{(0)}(X_t)+ \si_t(X_t)b_t^{(1)}(X_t, P_t^*\mu)\big\}\d t+\si_t(X_t)\d W_t,\ \ t\in [0,T].\end{equation}
Then $(P_t^\mu)^*\gg=\L_{X_t}.$ Let
$$\xi_s^\vv:=  b_s^{(1)}(X_s, P_s^* \pi_\vv)- b_s^{(1)}(X_s, P_s^*\mu),\ \ \ s\in [0,T].$$
 By combining \eqref{CP0} with \eqref{D**} and \eqref{DC1**}, we find a constant $c_1>0$ such that
\beq\label{*1} |\xi_s^\vv|\le K_s  \|P_s^* \pi_\vv-P_s^* \mu\|_{k,var}
\le c_1K_s(\mu+\nu)(|\cdot|^k).\end{equation}
Then by Girsanov's theorem,
\beq\label{W1} W_t^\vv:= W_t-\int_0^t  \xi_s^\vv \d s,\ \ t\in [0,T]\end{equation} is a Brownian motion under the weighted probability $R_T^\vv\d\P$, where
\beq\label{W2} R_t^\vv:= \e^{\int_0^t\<\xi_s^\vv,\d W_s\>-\ff 1 2 \int_0^t|\xi_s^\vv|^2\d s}, ~~t\in[0,T].\end{equation}
Reformulating the SDE \eqref{E01} as
$$\d X_t=  \big\{b_t^{(0)}(X_t)+ \si_t(X_t)b_t^{(1)}(X_t, P_t^*\pi_\vv)\big\}\d t+\si_t(X_t)\d W_t^\vv,\ \ t\in [0,T],$$ by the weak uniqueness and \eqref{DC1*}, we find some constant $c_2>0$ such that
\beq\label{LB1} \beg{split} &\|(P_t^{\pi_\vv})^*\gg- (P_t^\mu)^*\gg\|_{k,var}=\sup_{|f|\le 1+|\cdot|^k}
\big|\E[f(X_t)(R_t^\vv-1)]\big|\\
&\le \E[(1+|X_t|^k)|R_t^\vv-1|]\le \E\Big[\ss{\E((1+|X_t|^k)^2|\F_0)\E(|R_t^\vv|^2-1|\F_0)}\Big]\\
&\le c_2 (1+\E|X_0|^k) \Big\| \ss{\E(|R_t^\vv|^2-1|\F_0)}\Big\|_{L^\infty(\P)}\\
&=c_2 (1+\gg(|\cdot|^k)) \Big\| \ss{\E(|R_t^\vv|^2-1|\F_0)}\Big\|_{L^\infty(\P)}.\end{split}\end{equation}
By \eqref{*1} and that $\e^r-1\le r \e^r$ for $r\ge 0$, we find constant $c_3>0$ such that
\beg{align*}&\E(|R_t^\vv|^2-1|\F_0)\le \e^{c_3 \int_0^t \|P_s^*\pi_\vv-P_s^*\mu\|_{k,var}^2\d s}-1\\
&\le c_3 \e^{c_3 \mu(|\cdot|^k)^2+c_3 \nu(|\cdot|^k)^2}\int_0^t K_s^2\|P_s^*\pi_\vv-P_s^*\mu\|_{k,var}^2\d s.\end{align*}
Noting that \eqref{D**} implies
\beq\label{1**} P_s^*\pi_\vv-P_t^*\mu= (P_s^{\pi_\vv})^*\pi_\vv-P_t^*\mu=   (P_s^{\pi_\vv})^*\mu- P_t^*\mu +\vv  (P_t^{\pi_\vv})^* (\nu-\mu),\end{equation}
by the above estimate and \eqref{DC1**},
we find a constant $c_1(\mu,\nu)>0$ increasing in $(\mu+\nu)(|\cdot|^k) $ such that
\beq\label{*2}  \E(|R_t^\vv|^2-1|\F_0)\le c_1(\mu,\nu) \int_0^t K_s^2\|(P_s^{\pi_\vv})^*\mu- P_s^*\mu\|_{k,var}^2\d s
+c_1(\mu,\nu)\vv^2,\ \ t\in [0,T].\end{equation}
Combining this with \eqref{LB1} for $\gg=\mu$, and applying Gronwall's lemma, we find a constant $c_2(\mu,\nu)>0$ increasing in $(\mu+\nu)(|\cdot|^k) $  such that
$$\|(P_t^{\pi_\vv})^*\mu- (P_t^\mu)^*\mu\|_{k,var}^2\le c_2(\mu,\nu) \vv^2,\ \ t\in [0,T], \vv\in (0,1).$$
This together with \eqref{LB1} finishes the proof.

\end{proof}

\beg{lem}\label{L2} Assume  {\bf (H)},   let $\mu,\nu\in \scr P_k$ and $ \pi_\vv=(1-\vv)\mu+\vv\nu, \vv\in (0,1)$. Then
$$ \eta_t^{\vv,\mu,\nu}:= \ff 1 \vv \big(b_t^{(1)}(X_t^\mu, P_t^*\pi_\vv)-b_t^{(1)}(X_t^\mu,P_t^*\mu)\big)$$
satisfies
\beq\label{*3} \beg{split} &\eta_t^{\vv,\mu,\nu} = P_t^\mu\big(\tt D^E b_t^{(1)}(X_t^\mu, P_t^*\mu)\big)(\nu)- P_t\big(\tt D^E b_t^{(1)}(X_t^\mu, P_t^*\mu)\big)(\mu) \\
&\quad + \E\bigg[\tt D^E b_t^{(1)}(z, P_t^*\mu)(X_t^\mu)
\int_0^t \big\<\eta_s^{\vv,\mu,\nu}, \d W_s\big\>\bigg]_{z=X_t^\mu}+o_t(\vv),\ \ t\in [0,T]\end{split}\end{equation}
where  $\{o_t(\vv)\}_{t\in [0,T]}$ is   progressively measurable such that
\beq\label{04}  \lim_{\vv\downarrow 0} \sup_{t\in [0,T]}\ff{\|o_t(\vv)\|_\infty}{K_t}   =0.\end{equation}
\end{lem}
\beg{proof} (a)  By \eqref{D**} we have
\beq\label{**0}P_t^*\pi_\vv=(P_t^{\pi_\vv})^*\pi_\vv= (P_t^{\pi_\vv})^*\mu+\vv(P_t^{\pi_\vv})^*(\nu-\mu),\ \ t\in [0,T], \vv\in (0,1).\end{equation}
Next, by \eqref{LF} for $f=b_t(x,\cdot)$, we obtain
\beg{align*} &\eta_t^{\vv,\mu,\nu}=\ff 1 \vv \int_0^1\d r\int_{\R^d}\tt D^E b_t^{(1)}(X_t^\mu, r P_t^*\pi_\vv+(1-r)P_t^*\mu)\d (P_t^*\pi_\vv-P_t^*\mu)\\
&= \int_0^1\d r \int_{\R^d} \tt D^E b_t^{(1)}(X_t^\mu, r P_t^*\pi_\vv+(1-r)P_t^*\mu)\d \big((P_t^{\pi_\vv})^*(\nu -\mu)\big)\\
&\quad + \ff 1 \vv \int_0^1 \d r \int_{\R^d} \big[\tt D^E b_t^{(1)}(X_t^\mu, r P_t^*\pi_\vv+(1-r)P_t^*\mu)\\
&\qquad\qquad \qquad \qquad - \tt D^E b_t^{(1)}(X_t^\mu, r P_t^*\mu)\big]\d\big((P_t^{\pi_\vv})^*\mu-(P_t^\mu)^*\mu\big)\\
&\quad +  \ff 1\vv  \int_{\R^d} \tt D^E b_t^{(1)}(X_t^\mu,  P_t^*\mu)\d \big((P_t^{\pi_\vv})^*\mu -
(P_t^\mu)^*\mu\big).\end{align*}
Thus, denoting the right hand side of \eqref{*3} by $\tt \eta_t^{\vv,\mu,\nu}$, we obtain
$$\eta_t^{\vv,\mu,\nu}-\tt\eta_t^{\vv,\mu,\nu}= o_t(\vv),$$
where
\beg{align*}&o_t(\vv):= I_1(\vv,t)+ I_2(\vv,t)+I_3(\vv,t)+I_4(\vv,t),\\
& I_1(\vv,t):=\int_0^1\d r \int_{\R^d} \big[\tt D^E b_t^{(1)}(X_t^\mu, r P_t^*\pi_\vv+(1-r)P_t^*\mu)\\
&\qquad\qquad \qquad \qquad \qquad - \tt D^E b_t^{(1)}(X_t^\mu, P_t^*\mu)\big]\d \big((P_t^{\mu})^*(\nu-\mu)\big)\\
&I_2(\vv,t):= \int_0^1\d r \int_{\R^d}  \tt D^E b_t^{(1)}(X_t^\mu, P_t^*\mu) \d \big((P_t^{\pi_\vv})^*(\nu-\mu)-
(P_t^{\mu})^*(\nu -\mu)\big)\\
&I_3(\vv,t):= \ff 1 \vv \int_0^1 \d r \int_{\R^d} \big[\tt D^E b_t^{(1)}(X_t^\mu, r P_t^*\pi_\vv+(1-r)P_t^*\mu)\\
&\qquad\qquad \qquad \qquad \qquad \qquad - \tt D^E b_t^{(1)}(X_t^\mu,  P_t^*\mu)\big] \d\big((P_t^{\pi_\vv})^*\mu-(P_t^\mu)^*\mu\big)\\
&I_4(\vv,t):= \ff 1 \vv \int_{\R^d} \tt D^E b_t^{(1)}(X_t^\mu,  P_t^*\mu)\d \big((P_t^{\pi_\vv})^*\mu -
(P_t^\mu)^*\mu\big)\\
&\qquad \qquad \qquad - \E\bigg[\tt D^E b_t^{(1)}(z, P_t^*\mu)(X_t^\mu)
\int_0^t \big\<\eta_s^{\vv,\mu,\nu}, \d W_s\big\>\bigg]_{z=X_t^\mu}.\end{align*}
So, it suffices to verify \eqref{04}.
In the following, $c(\mu,\nu)$ stands for a constant depending on $\mu$ and $\nu$, which might be different in different places.

(b) By \eqref{DC3} and \eqref{DC1**},  we have
\beg{align*}& \big|\tt D^E b_t^{(1)}(X_t^\mu, r P_t^*\pi_\vv+(1-r)P_t^*\mu)
- \tt D^E b_t^{(1)}(X_t^\mu, P_t^*\mu)\big|\\
&\le K_t\aa\big(\|P_t^*\pi_\vv-P_t^*\mu\|_{k,var}\big)\big(1+|\cdot|^k
+\mu(|\cdot|^k)+\nu(|\cdot|^k)\big).\end{align*}
Next, Lemma \ref{L1} and \eqref{**0} imply
\beq\label{**a}\|P_t^*\pi_\vv-P_t^*\mu\|_{k,var}\le c(\mu,\nu)\vv,\ \ t\in [0,T],\vv\in (0,1).\end{equation}
Thus,
\beq\label{**1} \beg{split}&\big|\tt D^E b_t^{(1)}(X_t^\mu, r P_t^*\pi_\vv+(1-r)P_t^*\mu)
- \tt D^E b_t^{(1)}(X_t^\mu, P_t^*\mu)\big|\\
&\le K_t\aa\big(c(\mu,\nu) \vv\big) \big(1+|\cdot|^k+\mu(|\cdot|^k)+\nu(|\cdot|^k)\big),\ \ t\in [0,T],\vv\in (0,1).\end{split}\end{equation}
This  together with    \eqref{DC1**} and Lemma \ref{L1} implies
$$ |I_1(\vv,t)|+|I_3(\vv,t)|\le c(\mu,\nu) K_t \aa\big(c(\mu,\nu) \vv\big),\ \ \vv\in (0,1), t\in [0,T].$$
Moreover,  \eqref{DC2}  and Lemma \ref{L1} imply
\beg{align*} |I_2(\vv,t)|\le K_t\big\|(P_t^{\pi_\vv})^*(\nu-\mu)-(P_t^{\mu})^*(\nu -\mu)\big\|_{k,var}
\le c(\mu,\nu)K_t \vv,\ \ \vv\in (0,1), t\in [0,T].\end{align*}
  Therefore,
\beq\label{I12} \lim_{\vv\downarrow 0}
\sum_{i=1}^3\sup_{t\in [0,T]} \ff{\|I_i(\vv,t)\|_\infty}{K_t+1}   =0.\end{equation}

(c) To estimate $I_4(\vv,t)$, we use Girsanov's transform as in the proof of Lemma \ref{L1}. Note that $X_t^\mu$ solves \eqref{E01} for $X_0=X_0^\mu$, and the SDE \eqref{E01} can be  reformulated as
$$\d X_t^\mu= \big\{b_t^{(0)}(X_t^\mu)+\si_t(X_t^\mu)b_t^{(1)}(X_t^\mu, P_t^*\pi_\vv)\big\}\d t+\si_t(X_t^\mu)\d W_t^\vv,$$
where $W_t^\vv$ is a Brownian motion under the weighted probability $R_T^\vv\d\P$ as defined in \eqref{W1} and \eqref{W2}, so that
\beq\label{XI} R_t^\vv=\e^{\int_0^t\<\xi_s^\vv,\d W_s\>-\ff 1 2\int_0^t |\xi_s^\vv|^2\d s},\ \  \xi_s^\vv= \vv \eta_s^{\vv,\mu,\nu}. \end{equation}Then by Girsanov's theorem,
we obtain
\beq\label{A4} |I_4(\vv,t)| =\inf_{c\in\R}\bigg|\E\bigg[\big(\tt D^E b_t^{(1)}(z, P_t^*\mu)(X_t^\mu)-c\big)\Big(\ff{R_t^\vv-1}\vv-\int_0^t
\big\<\eta_s^{\vv,\mu,\nu},\d W_s\big\>\Big)\bigg]_{z=X_t^\mu}\bigg|.\end{equation}
Noting that $|\e^r-1-r|\le \ff 1 2(1+\e^r)r^2$ holds for $ r\in \R,$ we deduce form \eqref{XI} that
\beq\label{A5}\beg{split} & \bigg|\ff{R_t^\vv-1}\vv-\int_0^t
\big\<\eta_s^{\vv,\mu,\nu},\d W_s\big\>\bigg|\\
&\le \ff\vv  2\big(R_t^\vv+1\big)\bigg| \int_0^t\<\eta_s^{\vv,\mu,\nu},\d W_s\>-
\ff \vv 2 \int_0^t |\eta_s^{\vv,\mu,\nu}|^2\d s\bigg|^2.\end{split}\end{equation}
Next, \eqref{CP0} and Lemma \ref{L1} yield
\beq\label{A55}\int_0^T\sup_{\vv\in (0,1)} \|\eta_s^{\vv,\mu,\nu}\|_\infty^2\d s \le c(\mu,\nu).\end{equation}
 By combining this with \eqref{DC1*}, \eqref{DC1**}, \eqref{A4}, \eqref{A5},     and H\"older's inequality,
 we derive
 \beg{align*} & \|I_4(\vv,t)\|_\infty
\le \vv K \E\bigg[\big(1+|X_t^\mu|^k\big)\big(R_t^\vv+1\big)\bigg|\int_0^t \big\<\eta_s^{\vv,\mu,\nu},\d W_s\big\>-
\ff \vv 2 \int_0^t  |\eta_s^{\vv,\mu,\nu}|^2\d s\bigg|^2\bigg]\\
&\le K_t\vv\E\bigg\{ \big(\E[(1+|X_t^\mu|^{2k})|\F_0]\big)^{\ff 12} \big(\E[(1+R_t^\vv)^4|\F_0]\big)^{\ff {1} {4}}\\
&\qquad\qquad \qquad \times
  \bigg(\E\bigg[\bigg|\int_0^t \big\<\eta_s^{\vv,\mu,\nu},\d W_s\big\>-
\ff \vv 2 \int_0^t  |\eta_s^{\vv,\mu,\nu}|^2\d s\bigg|^{4}\bigg|\F_0\bigg]\bigg)^{\ff 1 4}\bigg\}\\
&\le c(\mu,\nu)K_t\vv.\end{align*}
 Combining this with \eqref{I12} we verify \eqref{04}, and hence finish the proof.

\end{proof}

\section{Proof of Theorem \ref{T1}}

\beg{proof}[Proof of Theorem \ref{T1}(1)]
 To prove  Theorem \ref{T1}(1), let
$$\scr M:=\bigg\{\eta: [0,T]\times \OO\to \R^d \ \text{is progressively\ measurable},  \int_0^T K_s^2\|\eta_s\|_\infty^2\d s <\infty\bigg\}.$$
It is a complete metric space under the distance
$$\rr_\ll(\eta,\tt\eta):= \bigg(\int_0^T\e^{-\ll t}\|\eta_t-\tilde{\eta_t}\|_\infty^2\d t\bigg)^{\ff 1 2}$$
for each constant $\ll>0.$ By the fixed point theorem, it suffices verify the $\rr_\ll$-contraction for some $\ll>0$ of the map
$\Phi:\ \scr M\to\scr M$ defined by
\beg{align*}\Phi_t(\eta):= &\,P_t^\mu\big\{\tt D^E b_t^{(1)}(X_t^\mu, P_t^*\mu)\big\}(\nu)-  P_t\big(\tt D^E b_t^{(1)}(X_t^\mu, P_t^*\mu)\big)(\mu)\\
&\quad +\E\bigg[\tt D^E b_t^{(1)}(z, P_t^*\mu)(X_t^\mu)\int_0^t \big\<\eta_s, \d W_s\big\>\bigg]_{z=X_t^\mu}.\end{align*}
Let $\eta\in \scr M$.  We first prove that $\Phi(\eta)\in \scr M.$
By \eqref{DC2} and \eqref{DC1**},   we find a constant $c_1>0$ such that
\beq\label{*5}\beg{split}&|\Phi_t(\eta)|\le   \big|P_t^\mu\big\{\tt D^E b_t^{(1)}(X_t^\mu, P_t^*\mu)\big\}(\nu) - P_t\big(\tt D^E b_t^{(1)}(X_t^\mu, P_t^*\mu)\big)(\mu) \big|\\
&\quad\qquad \qquad  +  \bigg|\E\bigg[\tt D^E b_t^{(1)}(z, P_t^*\mu)(X_t^\mu)\int_0^t \big\<\eta_s, \d W_s\big\>\bigg]\bigg|\\
&\le c_1 K_t(1+(\mu+\nu )(|\cdot|^k)\big) \\
&\qquad  + K_t\E\bigg(\Big(\E \big[1+ |X_t^\mu|^{2k}\big|\F_0\big]\Big)^{\ff 12} \bigg(\E\bigg|\int_0^tK_s^2
  |\eta_s|^2 \d s\bigg|^{\ff {k}{2(k-\theta)}}\bigg|\F_0\bigg)^{\ff 12}\bigg)\\
&\le  c_1 K_t(1+(\mu+\nu )(|\cdot|^k)\big) + c_1\big(1+\mu(|\cdot|^k)\big) \bigg(\int_0^t
K_s^2\|\eta_s\|_\infty^2\d s\bigg)^{\ff 1 2},\ \ t\in [0,T].\end{split}\end{equation} Then $\Phi(\eta)\in \scr M.$
Similarly, for any $\eta,\tt\eta\in \scr M$,  we have
\beg{align*} &|\Phi_t(\eta)-\Phi_t(\tt\eta)|\le \sup_{z\in\R^d} \E\bigg|\tt D^E b_t^{(1)}(z, P_t^*\mu)(X_t^\mu)\int_0^t \big\<(\eta_s-\tt\eta_s), \d W_s\big\>\bigg|\\
&\le c_1 K_t\big(1+\mu(|\cdot|^{k})\big)^{\ff\theta k}  \bigg(\int_0^t K_s^2\|\eta_s-\tt\eta_s\|^2_\infty\d s\bigg)^{\ff 1 2},\ \ t\in [0,T].\end{align*}
So, there exists a constant $c(\mu)>0$ such that for any $\ll>0$,
\beg{align*}&\rr_\ll(\eta,\tt\eta)^2 \le c(\mu)   \rr_\ll(\eta,\tt\eta)^2 \int_0^T K_t^2 \bigg(\int_0^t K_s^2 \e^{-\ll(t-s)}\d s\bigg)\d t.\end{align*}
 Therefore, when $\ll>0$ is large enough, $\Phi$ is contractive in $\rr_\ll$. The unique solution $\eta_t^{\mu,\nu}$ is the fixed point of $\Phi$, so that \eqref{*5} yields
 $$  \|\eta_t^{\mu,\nu}\|_\infty^2  \le 2c_1^2K_t^2 (1+(\mu+\nu) (|\cdot|^k)\big)^2+ 2c_1^2
 \big(1+\mu(|\cdot|^{k} )\big)^2   \int_0^tK_s^2
\|\eta_s^{\mu,\nu}\|_\infty^2 \d s,\ \ t\in [0,T],$$
hence the estimate \eqref{ET2} follows from Gronwall's lemma.\end{proof}

\beg{proof}[Proof of Theorem \ref{T1}(2)] Let $\mu,\nu\in \scr P_k$ and $f\in \D_q$ for some $q\in [0,k)$. By Theorem \ref{T1}(1) and Lemma \ref{L2}, it is easy to see that
\beq\label{A0} \lim_{\vv\downarrow 0} \sup_{t\in [0,T]} \ff{\|\eta_t^{\vv,\mu,\nu}-\eta_t^{\mu,\nu}\|_\infty}{K_t} =0.\end{equation}
Let $\pi_\vv=(1-\vv)\mu+\vv \nu, \vv\in (0,1)$ as before.  Then \eqref{**0} implies
\beq\label{A01}  \beg{split}&\ff{ P_tf(\pi_\vv)-P_tf(\mu)} \vv=\ff 1 \vv \int_{\R^d} f\d \big(P_t^*\pi_\vv- P_t^*\mu\big)\\
&=\ff 1 \vv\int_{\R^d} f\d\big((P_t^{\pi_\vv})^*\mu-(P_t^\mu)^*\mu\big)+ \int_{\R^d} f\d (P_t^{\pi_\vv})^*(\nu-\mu).\end{split}\end{equation}
By Lemma \ref{L1} we have
\beq\label{A02} \lim_{\vv\downarrow 0} \int_{\R^d} f\d (P_t^{\pi_\vv})^*(\nu-\mu)=
 \int_{\R^d} f\d (P_t^\mu)^*(\nu-\mu)= \int_{\R^d} P_t^\mu f \d (\nu-\mu). \end{equation}
To calculate the other term in \eqref{A01}, we make use of theorem with $R_t^\vv$ used in step $(c)$ in the proof of Lemma \ref{L2}. By   \eqref{DC1**} and \eqref{A55}, for $f\in \D_k$ we have
\beg{align*}&\sup_{\vv\in (0,1)}\E\bigg[|f(X_t^\mu)|\Big|\ff{R_t^\vv-1}\vv\Big|^2\bigg]\\
&\le \sup_{\vv\in (0,1)}\E\bigg[\Big(\E\big[|f(X_t^\mu)|\big|\F_0\big]\Big)^{\ff 1 2}\Big(\E\Big[\Big|\ff{R_t^\vv-1}\vv\Big|^4\Big|\F_0\Big]\Big)^{\ff 1 2}\bigg]<\infty.\end{align*}
So, by the dominated convergence theorem,
 \eqref{A5},   \eqref{A55} and  \eqref{A0} imply
\beg{align*} & \lim_{\vv\to\infty} \ff 1 \vv\int_{\R^d} f\d\big((P_t^{\pi_\vv})^*\mu-(P_t^\mu)^*\mu\big)\\
&= \lim_{\vv\to\infty}  \E\bigg[f(X_t^\mu)\ff{R_t^\vv-1}\vv\bigg]
= \E\bigg[f(X_t^\mu)\int_0^t \big\<\eta_s^{\mu,\nu},\d W_s\big\>\bigg].\end{align*}
This together with \eqref{A01} and \eqref{A02} implies the desired formula \eqref{BS}.

Finally, \eqref{BS}, \eqref{DC1*}, \eqref{DC1**} and \eqref{ET2},
we find   constants $c_1,c_2>0$ such that
\beg{align*} &\big|\tt D^E_\nu P_tf(\mu)\big|\le \|f\|_{k,var} \int_{\R^d}| (1+|x|^k) \d (P_t^\mu)^*(\mu+\nu)
+\|f\|_{k,var}\E\bigg[(1+|X_t^\mu|^k) \bigg|\int_0^t\<\eta_s^{\mu,\nu},\d W_s\>\bigg|\bigg]\\
&\le c_1\|f\|_{k,var} \big(1+(\mu+\nu)(|\cdot|^k)\big) + c_1\|f\|_{k,var}\E\bigg[\Big(\E\big((1+|X_t^\mu|^k)^2\big|\F_0\big)\Big)^{\ff 1 2}
\bigg(\int_0^tK_s^2\|\eta_s^{\mu,\nu}\|^2_\infty\d s\bigg)^{\ff 1 2}\bigg]\\
&\le c_2\|f\|_{k,var} \big(1+(\mu+\nu)(|\cdot|^k)\big) + c_2\|f\|_{k,var}\big(1+\mu(|\cdot|^k)\big)
\big(1+(\mu+\nu)(|\cdot|^k)\big)\e^{c\mu(|\cdot|^k)^{\ff{2\theta}k}}\ss t.\end{align*}
Therefore,  the derivative estimate \eqref{DE} holds for some constant $c>0$.
\end{proof}

\section{Two specific models}

In this section we apply Theorem \ref{T1} to two specific models: 1) non-degenerate DDSDEs with singular drift, 2) degenerate DDSDEs with weak monotone conditions.

\subsection{Singular case}
To measure the singularity, we recall some functional spaces introduced in \cite{XXZZ}.  For any $p\ge 1$, $L^p(\R^d)$ is the class of   measurable  functions $f$ on $\R^d$ such that
 $$\|f\|_{L^p(\R^d)}:=\bigg(\int_{\R^d}|f(x)|^p\d x\bigg)^{\ff 1 p}<\infty.$$
  For any $z\in\R^d$ and $r>0$,   let $B(z,r):=\{x\in\R^d: |x-z|< r\}$ be the open ball centered at $z$ with radius $r$.
For any $p,q>1$ and $t_0<t_1$, let $\tt L_q^p$ denote the class of measurable functions $f$ on $[t_0,t_1]\times\R^d$ such that
$$\|f\|_{\tt L_q^p}:= \sup_{z\in \R^d}\bigg( \int_{0}^{T} \|1_{B(z,1)}f_t\|_{L^p(\R^d)}^q\d t\bigg)^{\ff 1 q}<\infty.$$
 We will  take $(p,q)$ from the class
$$\scr K:=\Big\{(p,q): p,q\in  (2,\infty),\  \ff d p+\ff 2 q<1\Big\}.$$

\beg{enumerate}\item[{\bf (A)}] Let $k\in [0,\infty)$.
 \item[$(1)$]  $a:= \si\si^*$ is invertible with $\|a\|_\infty+\|a^{-1}\|_\infty<\infty$, where $\si^*$ is the transposition of $\si$,   and
$$ \lim_{\vv\to 0} \sup_{|x-y|\le \vv, t\in [0,T]} \|a_t(x)-a_t(y)\|=0.$$
\item[$(2)$]   $|b^{(0)}|\in \tt L_{q_0}^{p_0}$ for some  $(p_0,q_0)\in \scr K$.  Moreover, $\si_t$ is weakly differentiable such that
$$ \|\nn \si\|\le \sum_{i=1}^l f_i $$  holds for some $l\in \mathbb N$ and $0\le f_i \in \tt L_{q_i}^{p_i}$ with $  (p_i,q_i) \in  \scr K, 1\le i\le l.$
  \item[$(3)$]   for any $\mu\in C([0,T];\scr P_{k})$,  $b^{(1)}_t(x,\mu_t)$   is  locally bounded in $(t,x)\in [0,T]\times\R^d$. Moreover,  there exist  constants  $\vv,K>0$   such that
\beg{align*} (1+|x|)^{k-2} \big[\<b_t^{(1)}(x,\mu),x\>+\vv |x|\cdot |b_t^{(1)}(x,\mu)| \big] \le K \big(1+|x|^k+\mu(|\cdot|^k)\big)  \ \ x\in \R^d,\mu\in \scr P_k. \end{align*}
\end{enumerate}

\beg{thm}\label{T2} Assume {\bf (A)} and $(H_2)$. Then $(H_1)$ is satisfied, so that the assertions in Theorem \ref{T1} hold.\end{thm}

\beg{proof} By \eqref{CP0} which follows from $(H_2)$, we have
$$\big|\si_t(x)b_t^{(1)}(x,\mu)-\si_t(x)b_t^{(1)}(x,\nu)\big|\le K\|\si\|_\infty\|\mu-\nu\|_{k,var}.$$
Combining this with (1)-(3) in {\bf (A)}, we may apply Theorem 2.1 in \cite{R23} for $V:=(1+|\cdot|^2)^{\ff k 2}$ and the proof for the third assertion therein  to derive the well-posedness of  \eqref{DC'} and the estimate \eqref{DC1}.  Therefore, $(H_1)$ holds.\end{proof}

\subsection{Degenerate case}

In this part, we intend to apply Theorem \ref{T1} to   the following distribution dependent stochastic Hamiltonian system
for $X_t=(X_t^{(1)}, X_t^{(2)})\in R^{d_1}\times \R^{d_2}:$
$$\beg{cases} \d X_t^{(1)}= Z_t^{(1)}(X_t)\d t,\\
\d X_t^{(1)}=  Z_t^{(2)}(X_t,\L_{X_t})\d t+ \tt\si_t(X_t)\d B_t,\end{cases}$$
where  for $d:=d_1+d_2$,
$$Z=(Z^{(1)},Z^{(2)}): [0,T]\times \R^{d_1+d_2}\times \scr P_k\to \R^d,\ \ \tt\si: [0,T]\times\R^d\to \R^{d_2}\otimes \R^{d_2}$$
are measurable, $\tt \si$ is invertible, and $B_t$ is the $d_2$-dimensional Brownian motion.
Let $\tt B_t$ be a $d_1$-dimensional Brownian motion independent of $B_t$. This model is covered by \eqref{E0} for
$$b_t^{(0)}:= (Z_t^{(1)}, 0),\ \ b_t^{(1)}= (0, \tt\si_t^{-1}Z_t^{(2)}),\ \ \si_t={\rm diag}(0,\tt\si_t),
\ \ W_t=(\tt B_t,B_t).$$
So, below we come back to \eqref{E0} and make the following assumption.

\beg{enumerate}\item[{\bf (B)}] Let $k\in [0,\infty)$. For any $\gg_\cdot\in C([0,T];\scr P_k)$,
  $b_t(x,\gg_t)$ and $\si_t(x)$ is continuous in $x \in \R^d$ and locally bounded (i.e. bounded on bounded subsets of $[0,T]\times \R^d$),  and there exist   constants $K>0, p>k$ and a map $\kk: C([0,T];\scr P_k)\to (0,\infty)$   such that for any $t\in [0,T], x,y\in \R^d$ and $\gg_\cdot\in C([0,T];\scr P_k)$,
\beq\label{LM1} \beg{split}& (1+|x|^2)^{\ff {p-2} 2} \Big[2\<x, b_t(x,\gg_t)\> +\|\si_t(x)\|_{HS}^2+\ff{(p-2)|\si^*_t(x)x|^2}{|1+|x|^2}\Big]
 \\
&\qquad \le K \big[1+|x|^{p} +\gg_t(|x|^k)^{\ff p k}\big], \end{split}\end{equation}
\beq\label{LM2}2\<x-y, b_t(x,\gg_t)-b_t(y,\gg_t)\> +\|\si_t(x)-\si_t(y)\|_{HS}^2\le \kk(\gg_\cdot)|x-y|^2\log(2+|x-y|^{-1}).\end{equation}
\end{enumerate}
It is well known that \eqref{LM2} implies the pathwise uniqueness of the SDE \eqref{DC}, which together with the weak existence ensured by the local boundedness and continuity in $x$, implies the    well-posed  of \eqref{DC}, see for instance  \cite{FZ}. Moreover, by It\^o's formula for $(1+|X_t^{\gg_\cdot,x}|^2)^{\ff p 2},$ \eqref{LM1} implies the estimate \eqref{DC1}.  So, the following result follows from Theorem \ref{T1}.

\beg{thm}\label{T2} Assume {\bf (B)} and $(H_2)$. Then $(H_1)$ is satisfied, so that the assertions in Theorem \ref{T1} hold.\end{thm}

\end{document}